\documentstyle[11pt,fleqn]{article}

\title{\bf{The Garnir Relations for Weyl Groups }}

\author{\parbox{7cm}{\centering{\bf Sait ~Hal\i c\i o{\u g}lu \\}
Department of Mathematics\\ Ankara University\\ 06100 Tando{\u g}an
Ankara\\ Turkey\\}
\date{}
\vspace{0.1cm}}
\renewcommand{\Box}{\rule{2mm}{2mm}}
\renewcommand{\baselinestretch}{1.15}
\begin{document}
\maketitle
\renewcommand{\baselinestretch}{0.10}

\noindent
{\footnotesize
{\bf Abstract.} The so-called Garnir relations are very important in the
representation theory of the symmetric groups. In this paper,
 we give a possible extension of Garnir relations for Weyl groups
 in general.}

\renewcommand{\baselinestretch}{1.15}

\section{Introduction}
There are well known constructions
of the irreducible representations and of the
irreducible modules, called Specht modules, for
the symmetric groups $S_{n}$ which are based on  combinatorial
concepts connected with Young tableaux and tabloids  (see, e.g. [4]).
 In [5]
Morris described a possible extension of this work to
Weyl groups in general.  An alternative and improved approach
was described by the present author and Morris [1]. Later,
Hal\i c\i o{\u g}lu [2] develop the theory further and show how a $K$-basis for
Specht modules can be constructed in terms of standard
tableaux and tabloids.
The so-called Garnir relations are very important in the representation theory of the symmetric groups. The main aim of this paper
 is to give a possible extension of Garnir relations for Weyl groups in general .

\section{Some General Results On Weyl Groups}
We now state some results on Weyl groups which are required later .
 Any unexplained notation may be found in J E Humphreys [3],
 James and Kerber [4], Hal\i c\i o{\u g}lu and Morris [1],
 Hal\i c\i o{\u g}lu [2].

Let $V$ be $l$-dimensional Euclidean space over the real field  {\bf R} equipped with a positive
 definite inner product ( , ) for $ \alpha \in V~,~\alpha \neq 0$ , let
$\tau_{\alpha}$ be the $reflection$ in the hyperplane orthogonal to
$\alpha$ , that is , $\tau_{\alpha}$ is the linear transformation on $V$ defined by
\[ \tau_{\alpha}~(~v~)=v~-~2 ~\frac{(~\alpha~,~v~)}{(~\alpha~,~\alpha~)}~ \alpha\]
\noindent
for all $ v \in V$ . Let  $\Phi$ be a root system in $V$ and $\pi$ a simple
 system in $\Phi$ with corresponding positive system $\Phi^{+}$ and negative
 system $\Phi^{-}$ . Then , the $Weyl~group$ of $\Phi$
is the finite reflection group $ {\cal W}={\cal W}(\Phi)$,
 which is generated by the $\tau_{\alpha}, \alpha \in \Phi $.

\noindent
{\bf 2.1} We now give some of the basic facts presented in [1].

\noindent
Let $\Psi$  be  a subsystem of $\Phi$ with simple
system $J \subset \Phi^{+}$ and Dynkin diagram $\Delta$ and $\Psi=
\displaystyle \bigcup_{i=1}^{r}\Psi_{i} $, where $\Psi_{i}$ are the
indecomposable components of $\Psi$, then let $J_{i}$ be a simple
system in $\Psi_{i}$ ($i=1,2,...,r$) and $J = \displaystyle
\bigcup_{i=1}^{r}J_{i}$. Let $\Psi^{\perp}$ be the largest subsystem
in $\Phi$ orthogonal to $\Psi$ and let $J^{\perp} \subset
\Phi ^{+}$ be the simple system of $\Psi^{\perp}$.

Let $\Psi^{'}$ be  a subsystem of $\Phi$ which is contained in $\Phi
\setminus \Psi$, with simple system $J^{'} \subset \Phi^{+}$ and
Dynkin diagram $\Delta^{'}$,  $\Psi^{'}=\displaystyle
\bigcup_{i=1}^{s}\Psi_{i}^{'} $, where $\Psi_{i}^{'}$ are the
indecomposable components of $\Psi^{'}$; let $J_{i}^{'}$ be a
simple system in $\Psi_{i}^{'}$ ($i=1,2,...,s$) and $J = \displaystyle
\bigcup_{i=1}^{s}J_{i}^{'}$.  Let $\Psi^{'^{\perp}}$ be the largest
subsystem in $\Phi$ orthogonal to $\Psi^{'}$ and let $J^{'^{\perp}}
\subset
\Phi ^{+}$ be the simple system of $\Psi^{'^{\perp}}$.
 Let $\bar{J}$ stand for the ordered set $\{J_{1},J_{2},...,J_{r};
J_{1}^{'},J_{2}^{'},...,J_{s}^{'}\}$, where in addition the elements
in each $J_{i}$ and $J_{i}^{'}$ are also ordered. Let ${\cal T}_{\Delta}=\{ w\bar{J} \mid w\in {\cal W} \}$. The pair
$\{J,J^{'}\}$ is called a $useful~system $ in $\Phi$ if
${\cal W}(J)\cap {\cal W}(J^{'})=<e>$ and ${\cal W}(J^{\perp}) \cap
{\cal W}(J^{'^{\perp}})=<e>$. The
elements of ${\cal T}_{\Delta}$ are called $\Delta -tableaux$, the
$J$ and $J^{'}$ are called the $rows$ and the $columns$ of
$\{J,J^{'}\}$ respectively.
 Two $\Delta$-tableaux $\bar{J}$ and $\bar{K}$ are
$row-equivalent$, written $\bar{J}~\sim~\bar{K}$, if there exists
$w \in {\cal W}(J)$ such that $\bar{K}=w~\bar{J}$.
 The equivalence class which contains the $\Delta$-tableau $\bar{J}$ is
denoted by
$\{\bar{J}\}$ and is called a $\Delta-tabloid$. Let $\tau_{\Delta}$ be
set of all $\Delta$-tabloids. Then $\tau_{\Delta}=
\{ \{~d\bar{J}~\} \mid d\in D_{\Psi} \}$ , where
$D_{\Psi}=\{~w \in {\cal W} \mid w~(j) \in \Phi^{+}~ for~
all~j \in J~\}$ is a distinguished set of coset representatives of ${\cal W}
(\Psi)$ in ${\cal W}$ .
The group ${\cal W}$ acts on $\tau_{\Delta}$ as
$\sigma~\{\overline{wJ}\} =\{\overline{\sigma w J}\}$
for~all~$\sigma \in {\cal W}$.
Let $K$ be arbitrary field, let $M^{\Delta}$ be the
$K$-space whose basis elements are the $\Delta$-tabloids.
 Extend the
action of ${\cal W}$ on $\tau_{\Delta}$ linearly on $M^{\Delta}$,
then $M^{\Delta}$ becomes a $K{\cal W}$-module. Let
\begin {eqnarray*}
\kappa_{J^{'}}=\sum_{\sigma~
\in~{\cal  W}(J^{'})}~s~(~\sigma~)~
\sigma ~~and~~e_{J,J^{'}}=\kappa_{J^{'}}~\{~\bar{J}~\}
\end{eqnarray*}
\noindent
where $s(\sigma)=(-1)^{l(\sigma)}$ is the sign function and $l(\sigma)$ is
the length of $\sigma$.  Then $e_{J,J^{'}}$
is called the generalized $\Delta- polytabloid $ associated with
$\bar{J}$. Let $S^{J,J^{'}}$ be the subspace of $M^{\Delta}$
generated by $ e_{wJ,wJ^{'}}~$ where $~w\in~{\cal W}$.  Then
$S^{J,J^{'}}$
is called a $generalized~Specht~module$.
 A useful system $\{J,J^{'}\}$ in $\Phi$ is called a $good~system$ if $d~\Psi
\cap \Psi^{'} = \emptyset$ for $d \in D_{\Psi}$ then
$\{~\overline{dJ}~\}$ appears with non-zero coefficient in
$e_{J,J^{'}}$. If  $\{J,J^{'}\}$ is a good system ,
 then it was also proved in [1] that $S^{J,J^{'}}$ is  irreducible.

\noindent
{\bf 2.2}  The following are proved in [2].

A good system $\{J,J^{'}\}$ is called
a $very~good~system$ in $\Phi$ if
$d \leq d^{'}$ for all $d \in D_{\Psi} \cap D_{\Psi^{'}}$,
$d^{'} \in D_{\Psi}$ such that $d^{'}= d \sigma \rho $ , where $\rho \in
{\cal W}(J)~,~ \sigma \in
{\cal W}(J^{'})$.

\noindent
If $\{J,J^{'}\}$ is a very good system in $\Phi$ ,
then $\{e_{dJ,dJ^{'}} \mid d\in D_{\Psi}\cap D_{\Psi^{'}} \}$ is
linearly independent  over K.
\newline
A very good system $\{J,J^{'}\}$  is called a $perfect~system$ in
$\Phi $ if $\{ e_{dJ,dJ^{'}} \mid d \in D_{\Psi} \cap D_{\Psi^{'}} \}$ is
  a  basis  for  $S^{J,J^{'}}$.

\noindent
{\bf 2.3} The proof given by M . H . Peel  [6]  can be modified to
prove that if  $U$ and $V$ are subgroups of ${\cal W}$ and $Y = UV$,
then
\begin {eqnarray*}
(\sum_{\sigma ~\in ~U} ~s(\sigma)\sigma)(\sum_{\sigma ~\in ~V} ~s(\sigma)
\sigma)=\mid U \cap V \mid (\sum_{\sigma ~\in ~Y} ~s(\sigma)\sigma)
\end{eqnarray*}

\section{The Garnir Relations}

Let $\Phi$ be a root system and let $\{J,J^{'}\}$ be a
very good system in $\Phi$. Let $\overline{wJ}$ be a given $\Delta$-tableau,
where $ w \in {\cal W}$.  We want to find elements of $M^{\Delta}$
which annihilate the given $\Delta$-polytabloid $e_{wJ,wJ^{'}}$. By ( 2.1 ) if $w = d\rho$ , where $d \in D_{\Psi^{'}}$ and $\rho \in
{\cal W}(J^{'})$  we have $e_{wJ,wJ^{'}} ~=~s(\rho)~e_{dJ,dJ^{'}}$. Hence we may assume that $w \in D_{\Psi^{'}}$ .

\noindent
{\bf Lemma 3.1} $Let ~\{J,J^{'}\}~ be~ a ~  very~ good~ system~ in
~\Phi~ and~ let~ \Psi^{\ast}~ be~ a~ subsystem~of~ \Phi~ with$
 $simple~ system~ J^{\ast}.~ Let~d \in D_{\Psi^{'}}.~ Suppose
~there~ exists~ a~ 1:1~ and~onto~ function$~
\newline
$w \rightarrow w^{'}~ of~ the~ set~
 ~{\cal W}(J^{\ast}){\cal W}(d J^{'})~ with~ the~ following~ property:$
\newline
$( 3.1 )~ for~ each~  w \in {\cal W}(J^{\ast}){\cal W}(d J^{'}),~
there~ exists~ an~ element~ \rho_{w} \in {\cal W}(d J)~such~ that$~
\newline
$\rho_{w}^{2} = e ~,~s(\rho_{w}) = -1~,~w^{'} = w \rho_{w};~
further~(w^{'})^{'} = w.~ Then$
\begin {eqnarray*}
(\sum_{\sigma ~\in ~W(J^{\ast})} ~s(\sigma)\sigma)~e_{dJ,dJ^{'}}~=~0 ~.
\end{eqnarray*}

\noindent
{\bf Proof} By ( 2.3 )  we have
\begin {eqnarray*}
(\sum_{\sigma ~\in {\cal W}(J^{\ast})~} ~s(\sigma)\sigma)~e_{dJ,dJ^{'}}&=
&(\sum_{\sigma ~\in ~{\cal W}(J^{\ast})} ~s(\sigma)
\sigma)(\sum_{\sigma ~\in ~{\cal W}(dJ^{'})} ~s(\sigma)\sigma)\{\overline{dJ}\}\\
&=&\mid {\cal W}(J^{\ast}) \cap  {\cal W}(dJ^{'})\mid (\sum_{\sigma ~\in
{\cal W}(J^{\ast}){\cal W}(dJ^{'})~} ~s(\sigma)\sigma)\{\overline{dJ}\}
\end{eqnarray*}

But by ( 3.1 )  $(w^{'})^{'} = w $ , and thus the function is an involution.
If $w^{'} = w \rho_{w} = w$ then $\rho_{w} = e$ and  $s(\rho_{w}) = 1$ .
This is a contradiction.
Then the function is an involution without fixed points.
Suppose that $w_{1},w_{2} \in {\cal W}(J^{\ast}){\cal W}(d J^{'})$ and
$w_{1} \neq w_{2}$. If $w_{1}=w_{2}^{'}$ then $w_{1}^{'}=w_{2}$
(the function is an involution). Hence $w_{1}^{'} = w_{2} =
w_{1} \rho_{w_{1}} $ and $w_{2}^{'} = w_{1} = w_{2} \rho_{w_{2}}$. If $w_{1}= w_{1} \rho_{w_{1}}\rho_{w_{2}}$ then $\rho_{w_{1}}=
\rho_{w_{2}}$ and so $w_{1}=w_{2}$ . This is a contradiction.
But $\displaystyle{\bigcap_{i=1}^{2}\{~w_{i}~,~w_{i}^{'}~\}} =
\emptyset$ and ${\cal W}(J^{\ast}){\cal  W}(d J^{'})$ is a disjoint union of
2-elements sets $\{~w~,~w^{'}~\}$ .
Hence $\displaystyle{\sum_{\sigma ~\in {\cal W}(J^{\ast}){\cal W}(dJ^{'})~}
~s(\sigma)\sigma)\{\overline{dJ}\}}$ is a sum of terms of the
form $(s(\sigma)\sigma~+~s(\sigma^{'})\sigma^{'})\{\overline{dJ}\}$.

By ( 3.1 ) we have
\begin {eqnarray*}
(s(\sigma)\sigma~+~s(\sigma^{'})\sigma^{'})\{\overline{dJ}\}&=&(s(\sigma)\sigma~+~s(\sigma)s(\rho_{\sigma})\sigma \rho_{\sigma})\{\overline{dJ}\}\\
&=&(s(\sigma)\sigma~-~s(\sigma)\sigma \rho_{\sigma})\{\overline{dJ}\}\\
&=&(s(\sigma)\sigma)\{\overline{dJ}\}~-~(s(\sigma)\sigma \rho_{\sigma})\{\overline{dJ}\}\\
&=&(s(\sigma)\sigma)\{\overline{dJ}\}~-~(s(\sigma)\sigma )
\{\overline{dJ}\}\\
&=&0
\end{eqnarray*}

Then  \begin {eqnarray*}
(\sum_{\sigma ~\in ~{\cal W}(J^{\ast})} ~s(\sigma)\sigma)~e_{dJ,dJ^{'}}~=~0 ~. ~~~\Box
\end{eqnarray*}

If  $\sigma \in {\cal  W}(dJ^{'})$ then $\sigma e_{dJ,dJ^{'}} ~
=~s(\sigma)~e_{dJ,dJ^{'}}$ and so $s(\sigma)\sigma e_{dJ,dJ^{'}} ~=~
e_{dJ,dJ^{'}}.$
\noindent
Now let $H = \{~\sigma \in {\cal W}(J^{\ast}) ~\mid~s(\sigma)\sigma
e_{dJ,dJ^{'}}~=~e_{dJ,dJ^{'}}\} $ . Then $H = {\cal W}(J^{\ast}) \cap
{\cal W}(dJ^{'})$ and $H$ is a subgroup of ${\cal W}(J^{\ast})$.
Now let $C$ be a set of left coset representatives of $H$ in
${\cal W}(J^{\ast})$ such that $C$ contains the identity element .

\noindent
{\bf Theorem 3.2} $Let~ \{J,J^{'}\}~ be~ a ~ very~ good~ system~ in~ \Phi~
and~
d \in D_{\Psi^{'}}.~ Then$
\begin {eqnarray*}
e_{dJ,dJ^{'}}&=&-~(\sum_{\sigma ~\in C \setminus \{e\}} ~s(\sigma)\sigma)~e_{dJ,dJ^{'}}
\end{eqnarray*}

\noindent
{\bf Proof}
It follows directly from Lemma 3.1 by cancelling the factor
$\mid H \mid$ . $~~~~\Box$

\noindent
{\bf Definition 3.3} Let $\sigma_{1},\sigma_{2},...,\sigma_{k}$ be
coset representatives for H in ${\cal W}(J^{\ast})$ and let
\begin {eqnarray*}
G_{dJ,dJ^{'}}~=~\sum_{j=1}^{k} s(\sigma_{j}) \sigma_{j}
\end{eqnarray*}
$G_{dJ,dJ^{'}}$ is called $Garnir~element$ . $~~~~\Box$

\noindent
{\bf Example 3.4}
Let $\Phi={\bf G_{2}}$ with simple system
\[\pi = \{\alpha_{1},\alpha_{2}\}=\{\epsilon_{1}-\epsilon_{2},-2\epsilon_{1}+\epsilon_{2}+\epsilon_{3}\}\]
Let $a_{1}\alpha_{1}+a_{2}\alpha_{2}$ be denoted by
$a_{1}a_{2}$ and
$\tau_{\alpha_{1}},\tau_{\alpha_{2}}$
be denoted by $\tau_{1},\tau_{2}$ respectively.

Let $\Psi={\bf A_{1}+\tilde{A}_{1}}$ be the subsystem
of ${\bf G_{2}}$  with simple system $J=\{10,32\}$ and
 Dynkin diagram $\Delta$. Let $\Psi^{'}={\bf A_{1}}$  be another
subsystem of ${\bf G_{2}}$  with simple system $J^{'}=\{11\}$. Then
the $\Delta$-tabloids are:
\[\{\overline{J}\}=\{10,32;11\},\{\overline{\tau_{2}J}\}=\{11,31;10\} ,\{\overline{\tau_{1}\tau_{2}J}\}=\{21,01;-10\}\]
\noindent
Let $d=\tau_{1} \in D_{\Psi^{'}}$ . Then $\{\overline{dJ}\}=\{-10,32;21\}$. In this case
\newline
${\cal W}(dJ)=\{e,\tau_{1},\tau_{2}\tau_{1}\tau_{2}\tau_{1}\tau_{2},\tau_{2}\tau_{1}\tau_{2}\tau_{1}\tau_{2}\tau_{1} \}$ and ${\cal W}(dJ^{'}) = \{e,\tau_{1}\tau_{2}\tau_{1}\tau_{2}\tau_{1} \}$ .

If  $\Psi^{\ast}={\bf A_{2}}$ is the subsystem of ${\bf G_{2}}$
with simple system $J=\{10,21\}$, then
\newline
${\cal W}(J^{\ast)} = \{e,\tau_{1},\tau_{1}\tau_{2}\tau_{1}
\tau_{2}\tau_{1},\tau_{2}\tau_{1}\tau_{2},\tau_{1}\tau_{2}
\tau_{1}\tau_{2},\tau_{2}\tau_{1}\tau_{2}\tau_{1} \}$. Since
${\cal W}(J^{\ast}){\cal W}(d J^{'})= {\cal
W}(J^{\ast})$ and $w^{'} = w \tau_{1}$ and further $(w^{'})^{'}= w $
 for all $w \in {\cal W}(J^{\ast}){\cal W}(d J^{'})$, then Lemma
3.1 conditions are satisfied. Then
$H = {\cal W}(J^{\ast}) \cap {\cal W}(d J^{'})= {\cal W}(d J^{'}),
C = \{e,\tau_{1},\tau_{2}\tau_{1}\tau_{2} \}$ and the Garnir element
\newline
$G_{dJ,dJ^{'}} = e-\tau_{1}-\tau_{2}\tau_{1}\tau_{2}$ . By
Theorem 3.2 we have
\begin {eqnarray*}
e_{dJ,dJ^{'}}&=&\tau_{1}e_{dJ,dJ^{'}}+(\tau_{2}\tau_{1}\tau_{2})e_{dJ,dJ^{'}}\\
&=&e_{J,J^{'}}-e_{\tau_{2}J,\tau_{2}J^{'}}~~~~~~~\Box
\end{eqnarray*}

\end{document}